\documentclass{article}

\usepackage[preprint]{neurips_2026}

\usepackage{booktabs}
\usepackage{wrapfig}
\usepackage[utf8]{inputenc} 
\usepackage[T1]{fontenc}    
\usepackage{hyperref}       
\usepackage{url}            
\usepackage{booktabs}       
\usepackage{amsfonts}       
\usepackage{nicefrac}       
\usepackage{microtype}      
\usepackage{xcolor}         
\usepackage{graphicx}
\usepackage{subcaption}
\usepackage{bm}
\usepackage{graphicx}
\usepackage{amsmath}
\usepackage[version=4]{mhchem}
\usepackage{siunitx}
\usepackage{longtable,tabularx}
\usepackage[utf8]{inputenc}

\title{Accelerating the Simulation of Ordinary Differential Equations Through Physics-Preserving Neural Networks}

\author{%
    Andrew Tagg\thanks{PhD Candidate, Mechanical Engineering}, Andrew Frandsen\thanks{Undergraduate Research Assistant, Mechanical Engineering}, and Andrew Ning\thanks{Professor, Mechanical Engineering} \\
    Brigham Young University\\
    Provo, UT 84604 \\
    \texttt{\{atagg2, alfrand, aning\}@byu.edu}
}

\begin{document}

\maketitle

\begin{abstract}
    Numerical simulation of ordinary differential equations (ODEs) can be challenging when the system exhibits high accelerations and rapidly changing dynamics.  Under these conditions the ODE solver often needs to take very small time steps in order to resolve the solution accurately, resulting in increased computational cost.  In order to accelerate the simulation of these ODEs we present a novel methodology that uses a pseudo-invertible neural network to map system states into a high-dimensional latent-space.  The network is then trained so that the dynamics in this learned latent space are slow, and can be simulated with relatively few function calls.  Unlike existing neural methods, the latent dynamic equations are not learned from trajectory data, but derived from the original system equations and the chain rule.  This allows the method to generalize better than existing approaches because the derived equations are correct by construction.  In this work, we derive latent state equations of motion for any general ODE, and describe the loss function used to enforce slow time evolution of the latent states.  We then apply this technique to multiple example ODEs and show that these problems can be solved with $3$x to $20$x fewer function calls for the same accuracy when simulating in the learned latent space.  This reduction in cost could decrease computational demands for scientific simulations across engineering and physics applications.
\end{abstract}

\section{Introduction}

Ordinary differential equations (ODEs) are an essential mathematical tool for modeling how dynamic systems evolve.  The use of ODEs for analyzing and optimizing systems is widespread throughout all fields of physics and engineering.  Across all these applications, a common challenge when simulating ODEs occurs when a system exhibits rapidly changing dynamics.  When these high accelerations occur, the numerical solver must take very small steps in order to resolve the trajectories accurately.  This leads to increased computational cost, and simulations that can take hours or days to complete.  

The development of efficient numerical solution methods targeted toward fast changing dynamics is a well established field of research.  State of the art ODE solvers include Runge-Kutta methods \cite{kaps1979generalized, cooper1983additive}, Taylor methods \cite{tan2026efficient}, and multi-step methods \cite{butusov2020multistep}, all of which construct a high-order approximation of the local trajectory before taking a step.  For stiff systems, which have fast and slow dynamic modes, implicit methods are designed specifically to prevent numerical instability, but they often require more derivative evaluations to do so \cite{shampine1980implementation}.  All of these ODE solvers usually benefit from adaptive time stepping schemes \cite{lambda1998dynamical}.  These adaptive methods vary the time step of the solver based on an approximation of the local error of the solution.  In general, if the dynamics are smooth and slow changing, then the solver can safely take a larger time step.  If the dynamics start to change more rapidly, then the solver must shrink its time step to maintain accuracy.

Recently, deep learning has emerged as a viable approach for efficiently solving ODEs. Architectures like recurrent neural networks \cite{zhuang2025arecurrent}, LSTMs \cite{hong2025long}, and GRUs \cite{yuan2024application} naturally model sequential trajectory data without computing expensive derivatives. Neural ODEs \cite{tian2018neuralode} combine deep learning with advanced ODE solvers by representing the continuous time derivative with a neural network that can be simulated forward by any ODE solver. Other creative applications include learning mappings from coarse to fine discretizations \cite{cosso2026calorimeter}, improving accuracy for larger time steps \cite{mishra2018machinelearningframeworkdata}, and using autoencoders to map dynamic states to more desirable latent spaces \cite{qureshi2025autoencoders}.  Learning latent representations through autoencoders has led to improvements in dimensionality reduction \cite{guerin2026autoencoders}, Koopman linearization \cite{xiao2023deep}, neural ODE learning \cite{sholokhov2023physics}, and numerical stability \cite{erichson2019physicsinformedautoencoderslyapunovstablefluid}.

Neural networks can be a powerful tool for modeling ODEs efficiently, but their main limitation is the amount of data required to learn the underlying dynamics. In order to learn the dynamics of the system, the network must have access to a sufficiently large set of training trajectories. In many applications, obtaining enough data for the model to generalize is a challenge because the data is expensive to generate and/or the trajectory space is very high-dimensional.  This problem occurs in high-fidelity fluid simulations \cite{jinlong2017visualization}, astrodynamic simulations \cite{nejad2005acomparison}, financial markets \cite{falco2010algorithms}, and chemical reactions \cite{perini2014computationally}. In these domains where the dynamic equations are known but expensive to evaluate, computational cost is a major bottleneck. Advanced numerical solvers have led to drastic improvements, but speedups are still limited when the underlying dynamics are fast, and induce very small solver time steps.  Deep learning is a possible avenue to improve computational efficiency at run time, but the training process often requires a large data set that is prohibitively expensive to generate.

This paper presents a novel method to address these challenges and accelerate the simulation of ODEs.  The method employs a neural network architecture called the affine coupling layer \cite{dinh2016density} to define a pseudo-invertible mapping from the system's dynamic states to a higher-dimensional latent-space.  In this latent space, the dynamics are not learned from trajectory data, but derived directly from the original equations and the chain rule, making them automatically correct.  The goal of network training in this framework is to learn a mapping such that the derived latent dynamics evolve smoothly and slowly.  In this way, the simulation of the ODE is accelerated because the solver may safely take larger time steps.  Because the dynamic equations are correct by construction, this method can generalize better than the traditional approach of learning unknown dynamics from data.  

We first develop the methodology of this approach, deriving the equations of motion in the learned latent-space, and describe the loss function used to enforce slow latent dynamics.  Then we apply this technique to several benchmark examples, demonstrating accelerated ODE solutions with speedups ranging from $3$x to $20$x fewer function calls than the original simulations.  To evaluate the generalizability compared to existing methods, we also simulate these systems with a physics-informed neural ODE (PINODE) \cite{sholokhov2023physics} and show that while the PINODE performs well for trajectories similar to those provided during training, it struggles to extrapolate beyond the training data provided.  In contrast, the new method can predict trajectories outside ts training data with arbitrary accuracy. 

\section{Methods}


\subsection{Pseudo-Invertible Architecture}

Let $\mathbf{x} \in \mathbb{R}^n$ be the state vector of an ODE of interest.  We will learn a latent representation of the state-space by defining a mapping $\mathbf{z}(\mathbf{x})$, where $\mathbf{z} \in \mathbb{R}^{m > n}$.  As we will see, constructing $\mathbf{z}$ in a higher-dimensional space will allow more expressivity in the latent representation and improve the network's ability to learn dynamics that are slow in $\mathbf{z}$-space.   We desire this mapping to be pseudo-invertible, which is important not only for generalization of the mapping, but also for the validity of the derived equations of motion.  To achieve this, we will define the following mapping

\begin{equation}
    \label{eq:mapping}
    \mathbf{z} = \phi(\mathbf{A}\mathbf{x})
\end{equation}
where $\mathbf{A} \in \mathbb{R}^{m \times n}$ is a linear operation that lifts the original states to a higher dimension, and $\phi$ is a nonlinear invertible map defined by a neural network.  The pseudo-inverse of the mapping is

\begin{equation}
    \label{eq:inverse}
    \mathbf{x} = \mathbf{A}^{\dag}\phi^{-1}(\mathbf{z})
\end{equation}
where $\mathbf{A}^{\dag}$ is the Moore-Penrose pseudo-inverse defined by $\mathbf{A}^{\dag} = (\mathbf{A}^T \mathbf{A})^{-1}\mathbf{A}^T$.  For the function $\phi$, we use a neural network consisting of affine coupling layers \cite{dinh2016density}, which ensure that $\phi$ automatically has a well defined inverse $\phi^{-1}$. 

\subsection{Latent Dynamics}

\begin{figure}[hbt!]
\centering
\includegraphics[width=0.9\textwidth]{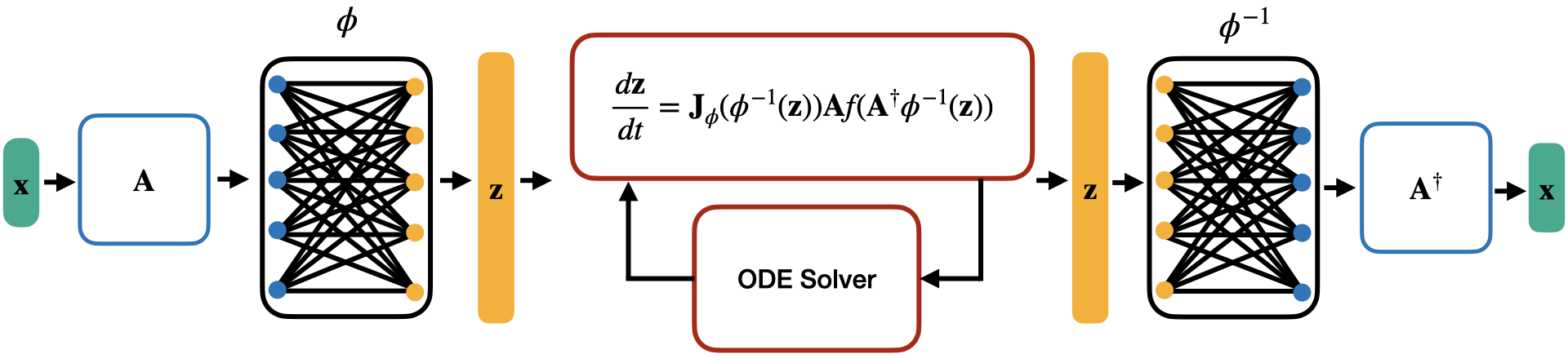}
    \caption{Procedure for simulating the dynamic system in the learned latent space, then reconstructing the result in the original space.}
\label{fig:architecture}
\end{figure}

Consider any dynamic system defined by a set of differential equations $\frac{d\mathbf{x}}{dt} = f(\mathbf{x})$.  Given these equations, and the mapping defined in Eq. \ref{eq:mapping}, we want to find an expression for $\frac{d\mathbf{z}}{dt}$.  We do this by applying the chain rule to $\mathbf{z}(\mathbf{x})$ 

\begin{equation}
    \label{eq:chain_rule}
    \frac{d \mathbf{z}}{dt} = \frac{d \mathbf{z}}{d \mathbf{x}} \frac{d \mathbf{x}}{dt}  
\end{equation}
We obtain $\frac{d \mathbf{z}}{d \mathbf{x}}$ by differentiating Eq. \ref{eq:mapping}. 

\begin{equation}
    \label{eq:f_sub}
    \frac{d \mathbf{z}}{dt} = \mathbf{J}_{\phi} \mathbf{A} \frac{d\mathbf{x}}{dt}
\end{equation}
where $\mathbf{J}_{\phi}$ is the Jacobian of $\phi$.  Both $\mathbf{J}_{\phi}$ and $\frac{d\mathbf{x}}{dt}$ are functions of $\mathbf{x}$.

\begin{equation}
    \label{eq:funcs_of_x}
    \frac{d \mathbf{z}}{dt}(\mathbf{x}) = \mathbf{J}_{\phi}(\mathbf{A}\mathbf{x}) \mathbf{A} f(\mathbf{x}) 
\end{equation}
Because we want the equations in terms of $\mathbf{z}$, we substitute in Eq. \ref{eq:inverse} to obtain

\begin{equation}
    \label{eq:z_dynamics}
    \frac{d \mathbf{z}}{dt}(\mathbf{z}) = \mathbf{J}_{\phi}(\phi^{-1}(\mathbf{z})) \mathbf{A} f(\mathbf{A}^{\dag}\phi^{-1}(\mathbf{z})) 
\end{equation}
These are the new equations of motion derived entirely in terms of $\mathbf{z}$. Note that instead of forming the full Jacobian $\mathbf{J}_{\phi}$, we compute a Jacobian-vector-product to evaluate the entire expression, which is more computationally efficient.  The full procedure for simulating the dynamic system in the new latent-space is depicted visually in Fig. \ref{fig:architecture}

Given an initial state in $\mathbf{x}$ space, we use $\mathbf{A}$ to lift the state into a higher-dimensional space, then pass it through the invertible neural network $\phi$. The resulting latent state $\mathbf{z}$ is simulated forward in time using Eq. \ref{eq:z_dynamics}. After simulation, we reconstruct the solution in $\mathbf{x}$-space by passing it through the inverse map $\phi^{-1}$ and projecting onto $\mathbf{x}$-space using the pseudo-inverse $\mathbf{A}^{\dag}$.

\subsection{Loss Function}

Now that we have an expression for the latent dynamics, we need to train $\phi$ and $\mathbf{A}$ to enforce slow time evolution of the latent states.  We do this by penalizing the $L2$ norm of the latent system's Jacobian matrix, which quantifies how fast the dynamics are changing in all directions.  Taking the average of the Jacobian norm across $N$ sample points, we have 

\begin{equation}
    \label{eq:spectral}
    L = \frac{1}{N} \sum_{i=1}^N \left| \left| \mathbf{J}_{\dot{\mathbf{z}}}(\mathbf{z}_i) \right| \right|_2^2
\end{equation}
where $\mathbf{J}_{\dot{\mathbf{z}}}(\mathbf{z}_i)$ is the Jacobian of $\frac{d \mathbf{z}}{dt}(\mathbf{z})$ evaluated at a sample point $\mathbf{z}_i$. In practice, forming the full Jacobian $J_{\dot{\mathbf{z}}}(\mathbf{z}_i)$ can be computationally expensive if the dimension of the latent-space is large.  Instead we use a surrogate for the norm by computing the average norm of $k$ random Jacobian-vector products.  Our final loss function is
\begin{equation}
    \label{eq:j_loss}
    L = \frac{1}{Nk} \sum_{i=1}^N \sum_{j = 1}^k ||\mathbf{J}_{\dot{\mathbf{z}}}(\mathbf{z}_i) \mathbf{v}_j||_2^2
\end{equation}
where $\mathbf{v}_i \in \mathbb{R}^m$ is a randomly generated vector.  In the following section, we explore several example problems, and apply this loss to the derived governing equations of the latent-space in order to smooth out the dynamics and ensure that trajectories evolve slowly in time.

\section{Results}

The following examples demonstrate the proposed approach across several dynamical systems. All experiments were performed on consumer-grade CPU hardware with 16 GB RAM.

\subsection{Linear System Example}
\label{sec:linear_system}

The first benchmark example was a simple $3$ state linear system of equations

\begin{equation}
    \label{eq:vortex_states}
    \begin{bmatrix} \dot{x}_1 \\ \dot{x}_2 \\ \dot{x}_3\end{bmatrix} = \begin{bmatrix} 33 & 17 & -70 \\ 42 & 18 & -80 \\ 37 & 18 & -75 \end{bmatrix}\begin{bmatrix} x_1 \\ x_2 \\ x_3 \end{bmatrix}
\end{equation}

This system has eigenvalues of $\bm{\lambda} = (-20, -2 + i, -2 - i)^T$.  Its trajectories exhibit an initial fast transient portion, followed by slower decay.  In order to slow down the fast dynamics, we trained a pseudo-invertible network to penalize the Jacobian loss in Eq. \ref{eq:j_loss} across $600$ random sample points.  

\begin{figure}[hbt!]
\centering
\begin{minipage}{0.49\textwidth}
\centering
\includegraphics[width=\textwidth]{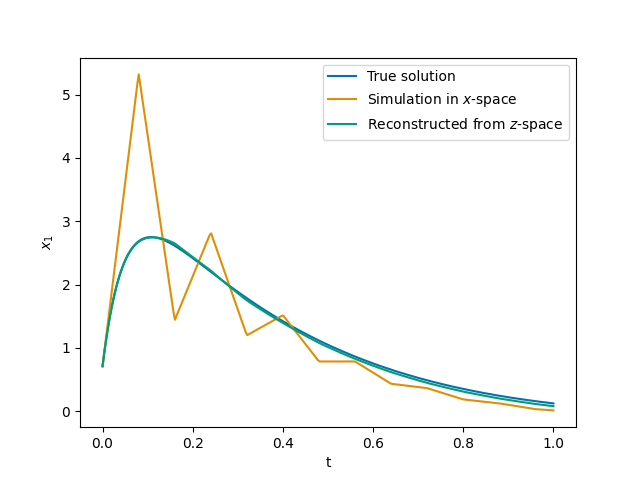}
\subcaption{Original state-space}
\end{minipage}\hfill
\begin{minipage}{0.5\textwidth}
\centering
\includegraphics[width=\textwidth]{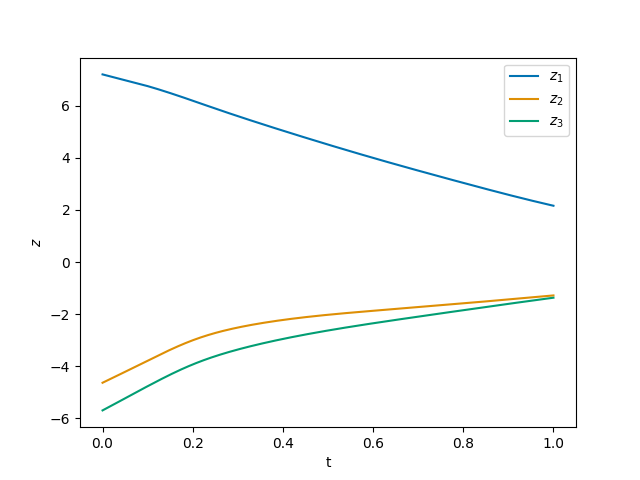}
\subcaption{Latent state-space}
\end{minipage}\hfill
    \caption{Linear system trajectories, plotted in the original state-space and the latent state space.  Trajectories are simulated with an Euler solver using a time step of $0.08$s.}
    \label{fig:linear_trajectory}
\end{figure}

After training, we simulated the system in both the original and learned latent-space using Euler integration with a time step of $0.08$s---significantly larger than what was required for the original simulation to be accurate. Fig. \ref{fig:linear_trajectory} shows the trajectories in $\mathbf{x}$ and $\mathbf{z}$ space.  For clarity, only one original state and three latent states were plotted, but the remaining states showed similar behavior. As shown, the fast dynamic response of the system is diminished in the latent space, making the $0.08$s Euler solution much more accurate in $\mathbf{z}$ than in $\mathbf{x}$.  While the original simulation overshot significantly, the latent simulation closely matched the true solution when reconstructed from $\mathbf{z}$ space.



\begin{figure}[hbt!]
\centering
\includegraphics[width=0.5\textwidth]{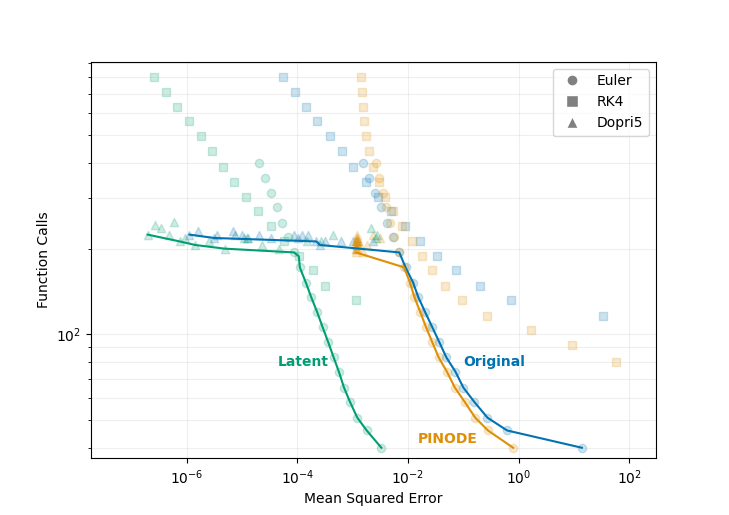}
    \caption{Function calls plotted against final mean squared error for the original, PINODE, and latent equations using Euler, RK4, and Dopri5 solvers.}
\label{fig:linear_work_precision}
\end{figure}

We repeated this experiment across a range of time steps and error tolerances for three solvers: Euler, RK4, and Dopri5 (adaptive). To compare with existing data-driven methods, we also trained a physics-informed neural ODE (PINODE) \cite{sholokhov2023physics} to learn the dynamics of the system.  We used the same $600$ sample points that were used to train the pseudo-invertible network as initial conditions to generate training trajectories for the PINODE. Fig. \ref{fig:linear_work_precision} compares latent, original, and PINODE simulations by plotting total function calls against final solution error on unseen test trajectories.  We plot data points for all solvers here, but the best solver performance for each error level is plotted as a continuous curve for better comparison. 


The speedups achieved by latent simulations vary with the desired precision. For moderate error tolerances (around $10^{-3}$), Euler latent simulations use about $5$x fewer function calls than original Dopri5 simulations. However, for tighter tolerances ($\le 10^{-6}$), Dopri5 performs similarly in both spaces. This suggests that for this linear system, the latent representation mainly enables simpler integration schemes to achieve competitive accuracy at moderate tolerances.  For tighter tolerances, Dopri5 already handles the short transient and long decay efficiently. As we will see in following examples, more complex nonlinear systems can benefit from latent representation even with Dopri5 at tighter tolerances.

The PINODE was also able to simulate the linear system and in some cases uses fewer function calls than the original simulations. However, there is a limit to the PINODE's accuracy on these test trajectories. Refining the solver step size reduces error down to about $10^{-3}$, but further refinement does not increase accuracy. This indicates that there is some irreducible error in the PINODE simulations.  More training data was likely required to generalize with higher accuracy.

\subsection{Vortex Particle Dynamics}
\label{sec:vortex_particle_dynamics}

\begin{wrapfigure}{r}{0.35\textwidth}
    \hyphenpenalty=10000
    \exhyphenpenalty=10000
    \vspace{-10pt}
    \centering
    \includegraphics[width=0.35\textwidth]{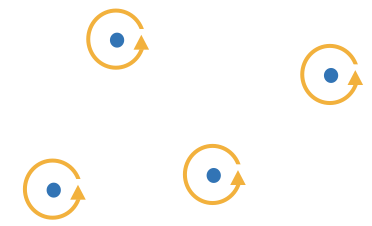}
    \caption{Depiction of the 2D vortex particle system.}
    \label{fig:vortices}
    \vspace{-10pt}
\end{wrapfigure}
Another ODE we considered was a system of 2D vortex particles. A vortex particle is a differential region of rotating fluid. As shown in Fig. \ref{fig:vortices}, each vortex induces a circulating velocity field surrounding the particle.  As these velocity fields combine, the particles move and interact with one another in complex ways.  The state vector for this system is a vector of all the positions of the particles $\mathbf{x} = [x_1, y_1, x_2, y_2, x_3, y_3, x_4, y_4]^T$.  The state derivatives are the velocities of each particle, which we can compute by summing the velocities induced by every other particle.  The velocity of the $i^{\text{th}}$ particle is 

\begin{equation}
    \label{eq:v_induced}
    \mathbf{V}_i = \sum_{j = 1}^4\frac{\mathbf{\Gamma}_j \times \mathbf{r}_{ij}}{2 \pi |\mathbf{r}_{ij}|}
\end{equation}
where $\mathbf{\Gamma}_j$ is the circulation of the $j^{\text{th}}$ vortex, and $\mathbf{r}_{ij} = [x_i - x_j, y_i - y_j]^T$.  The induced velocities are inversely proportional to the distance between the particles $|\mathbf{r}_{ij}|$.  This means that when the particles are far apart the dynamics are slow, but when they get close to each other the system is very fast. 

\begin{figure}[hbt!]
\centering
\begin{minipage}{0.5\textwidth}
\centering
\includegraphics[width=\textwidth]{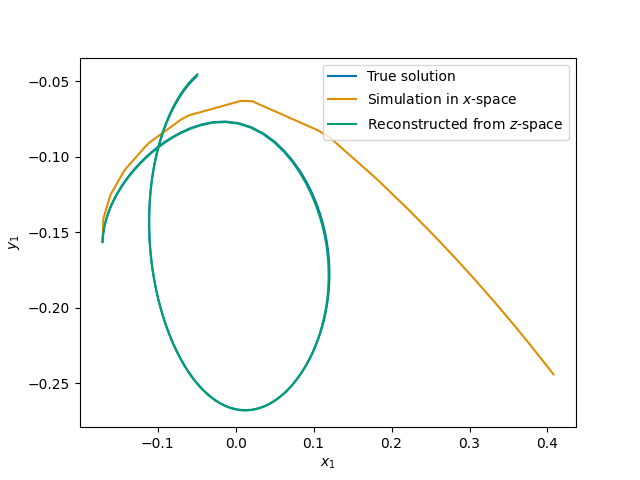}
\subcaption{Original state-space}
\end{minipage}\hfill
\begin{minipage}{0.5\textwidth}
\centering
\includegraphics[width=\textwidth]{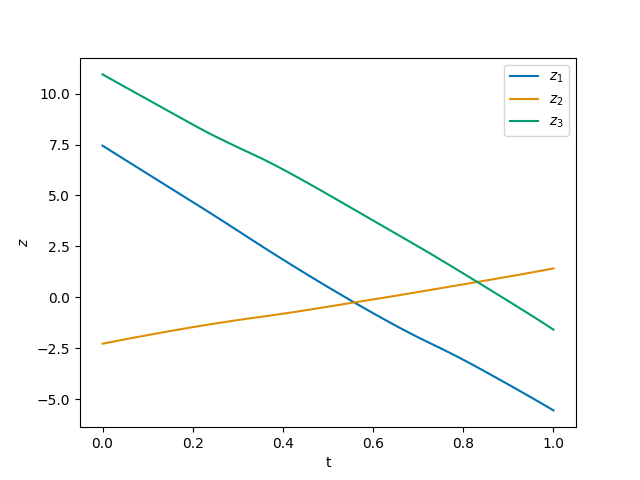}
\subcaption{Latent state-space}
\end{minipage}\hfill
    \caption{Vortex particle trajectory, plotted in the original state-space and the latent state space. Trajectories are simulated with an Euler solver using a time step of $0.06$s.}
    \label{fig:leapfrogging_trajectory}
\end{figure}

To reduce the number of function calls needed to simulate this system, we trained a pseudo-invertible network to smooth out the latent dynamics across $40$ training trajectories. We then simulated with an Euler time-step of $0.06$s.  The results of these simulations are shown in Fig. \ref{fig:leapfrogging_trajectory}, where the path of one particle is shown along with $3$ latent state trajectories. 

\begin{figure}[hbt!]
\centering
\begin{minipage}{0.5\textwidth}
\centering
\includegraphics[width=\textwidth]{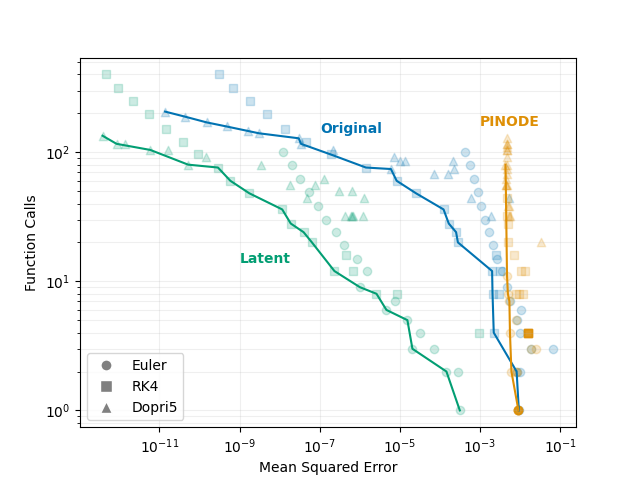}
\subcaption{Solver/method comparisons.}
\end{minipage}\hfill
\begin{minipage}{0.5\textwidth}
\centering
\includegraphics[width=\textwidth]{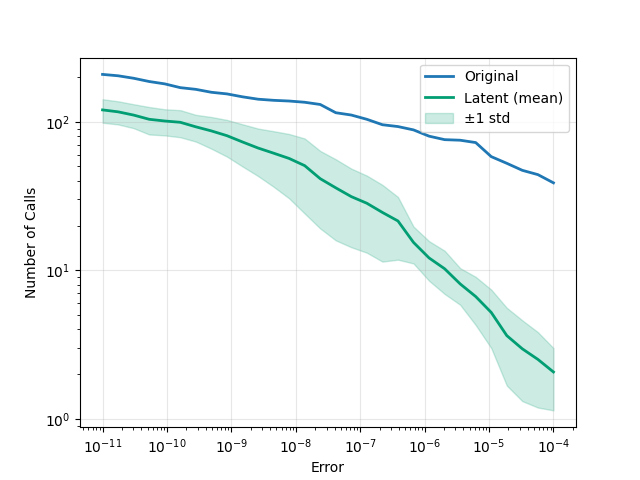}
\subcaption{Best solver results with $\pm 1$ standard deviation.}
\end{minipage}\hfill
    \caption{Function evaluations vs. mean squared error for the vortex particle simulations.}
    \label{fig:particle_work_precision}
\end{figure}

Notably, the curved particle trajectories are mapped to straight trajectories in latent-space, making the latent simulation much more accurate than the original. Fig. \ref{fig:particle_work_precision}(a) shows function evaluations vs. mean squared error for Euler, RK4, and Dopri5 solvers across original, latent, and PINODE simulations.  Fig. \ref{fig:particle_work_precision}(b) shows the best result at each error level, with $\pm 1$ standard deviation shown across $5$ training runs with random starting weights and random training samples.

%
%
The latent simulations consistantly use fewer function calls for the same accuracy as the original simulations. The largest speedups ($10$x to $20$x fewer function calls) occur for errors in the range of $10^{-5}$ to $10^{-3}$. For this example, the PINODE simulations produced errors consistently higher than $10^{-2}$, indicating it could not learn the correct dynamics from the sparse set of $40$ training trajectories.

\subsection{Vortex Lattice Model}
\label{sec:vortex_lattice_model}

\begin{wrapfigure}{r}{0.4\textwidth}
    \hyphenpenalty=10000
    \exhyphenpenalty=10000
    \vspace{-10pt}
    \centering
    \includegraphics[width=0.38\textwidth]{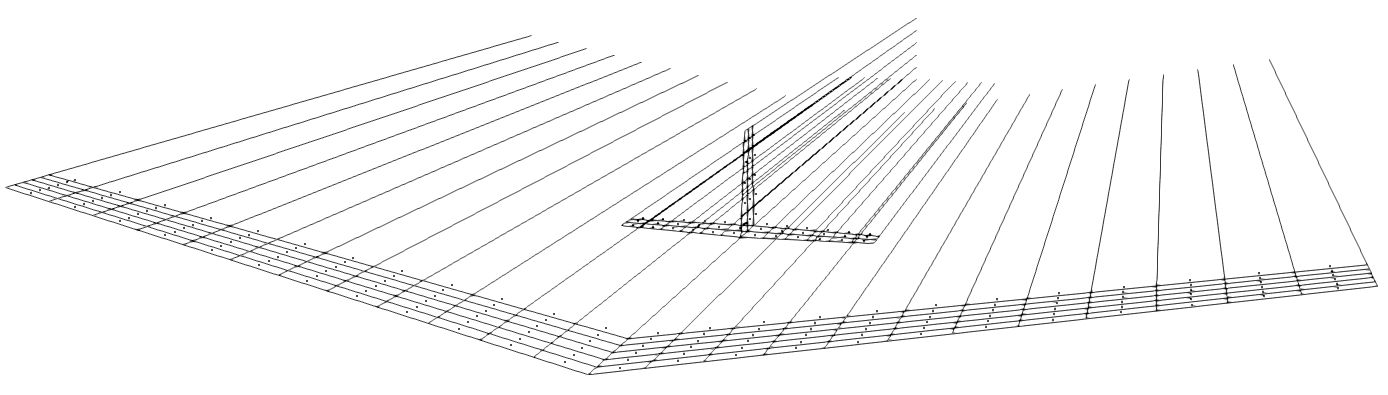}
    \caption{Depiction of the vortex lattice aerodynamic model.}
    \label{fig:wing_tail}
    \vspace{-10pt}
\end{wrapfigure}
The previous examples were relatively simple, and featured computationally inexpensive dynamic functions. This example demonstrates how the method scales to more complex problems using an aircraft modeled by a vortex lattice method (VLM). The VLM predicts aerodynamic forces by discretizing the wing surface into panels with vortex filaments.  Vortex filaments are similar to vortex particles, exept instead of a single point, we have a line of rotating fluid. Fig. \ref{fig:wing_tail} shows the VLM discretization with filaments extending behind the aircraft to model the wake.  Using these vortex filaments, we can compute the local velocities induced by each vortex filament on each panel of the aircraft model.  We then apply boundary conditions at each of these panels to ensure that the total local velocity is tangent to the surface everywhere.  These boundary conditions produce a linear system of equations that allows us to solve for the circulation strength $\mathbf{\Gamma}$ of each vortex filament.  Once this system is solved, we compute the aerodynamic force acting on each panel by applying the Kutta-Jouwkowski relation $\mathbf{F} = \rho \mathbf{V} \times \mathbf{\Gamma}$, where $\rho$ is the air density, and $\mathbf{V}$ is the local velocity at the panel. 

The forces acting on the aircraft can then be used to simulate the motion of the aircraft by applying classical rigid body mechanics.  In this example, we consider only 2D longitudinal motion, and the state vector consists of the horizontal and vertical velocity, pitch rate, and pitch angle of the aircraft. $\mathbf{x} = [v_x, \; v_z, \; \dot{\theta} \; \theta]^T$  We do not include the translational position of the aircraft as states in this formulation because they do not affect the dynamics of the aircraft, and can be obtained post-simulation by integrating the velocity signals.  With these states defined, the governing differential equations are
\begin{equation}
\label{eq:EOM}
    \frac{d\mathbf{x}}{dt} = f(\mathbf{x}) = \begin{bmatrix} 
                    \dot{v}_x \\
                    \dot{v}_z \\
                    \ddot{\theta} \\
                    \dot{\theta}
    \end{bmatrix} = 
    \begin{bmatrix} 
                    \frac{F_x(\mathbf{x})}{m} - \dot{\theta}v_z \\
                    \frac{F_z(\mathbf{x})}{m} - \dot{\theta}v_x \\
                    \frac{M_y(\mathbf{x})}{I_{yy}} \\
                    \dot{\theta}
    \end{bmatrix}
\end{equation}
where $F_x(\mathbf{x})$ and $F_z(\mathbf{x})$ are the net forces acting on the vehicle including VLM-computed aerodynamics, and gravity.  $M_y(\mathbf{x})$ is the pitching moment acting on the aircraft, which is calculated from the forces acting on each panel by summing moments about the center of gravity. 

\begin{figure}[hbt!]
\centering
\begin{minipage}{0.49\textwidth}
\centering
\includegraphics[width=\textwidth]{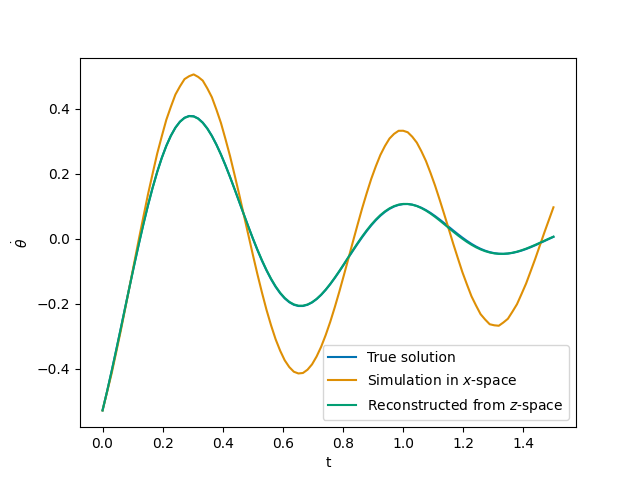}
\subcaption{Original state-space}
\end{minipage}\hfill
\begin{minipage}{0.5\textwidth}
\centering
\includegraphics[width=\textwidth]{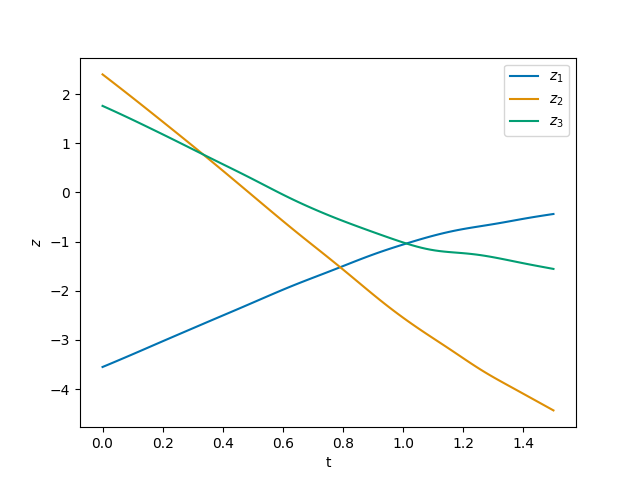}
\subcaption{Latent state-space}
\end{minipage}\hfill
    \caption{Original and Latent simulations for the VLM aircraft model.  Solutions are obtained using an Euler solver with a timestep of $0.03$s.}
    \label{fig:vlm_trajectory}
\end{figure}

For this example, we constructed a VLM model consisting of an $11$m wing, and a $2.75$m horizontal stabilizer.  We used $100$ panels to discretize the wing, and $25$ panels to discretize the tail.  We then trained a pseudo-invertible network to smooth out the dynamics of the aircraft model across $20$ sample trajectories.  Fig. \ref{fig:vlm_trajectory} shows the results of Euler simulations with a time step of $0.03$s in both the latent, and original space.  As shown, the Latent trajectories are effectively smoothed out, leading to higher accuracy when simulating in the latent-space. 

\begin{figure}[hbt!]
\centering
\includegraphics[width=0.5\textwidth]{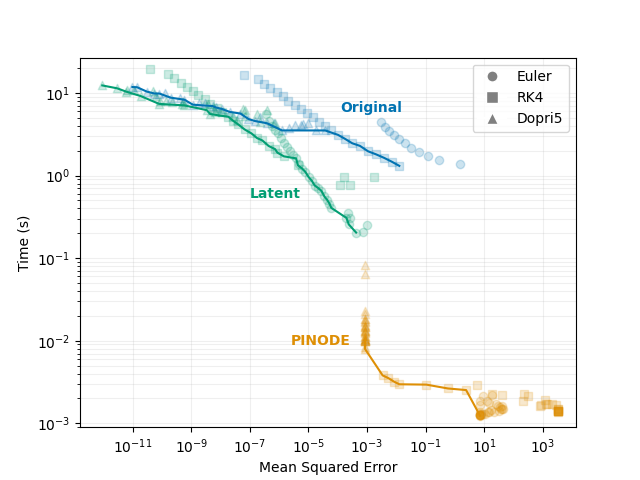}
    \caption{Function evaluations vs. mean squared error for the vortex lattice aircraft simulations using original, PINODE, and latent equations.  Solutions are obtained using Euler, RK4, and Dopri5 solvers.}
\label{fig:glider_speedups}
\end{figure}

In Fig. \ref{fig:glider_speedups}, we plot total simulation wall time vs. final error for the latent, original, and PINODE simulations using each ODE solver.  The speedups again vary depending on the desired accuracy.  For final errors in the range of $10^{-5}$ to $10^{-3}$, the original simulations take a few seconds to complete, whereas the latent simulations take a few tenths of a second, resulting in an overall speedup of $3$x to $10$x in this range of errors.  The PINODE simulations are $1$ to $2$ orders of magnitude faster than the original and latent simulations because the expensive function $f(\mathbf{x})$ is replaced by a relatively fast neural network.  This is a distinct advantage of the PINODE, however it could not achieve errors below $10^{-3}$ on these test trajectories. 

\subsection{Sensitivity Analyses}

In order to further understand the capabilities and limitations of the new method, we conducted $3$ studies which investigated the performance of the pesudo-invertible networks under various conditions. In these studies, we examined the effect of sample size, latent dimension, and problem complexity on the performance of the network. 

\begin{figure}[hbt!]
\centering
\begin{minipage}{0.49\textwidth}
\centering
\includegraphics[width=\textwidth]{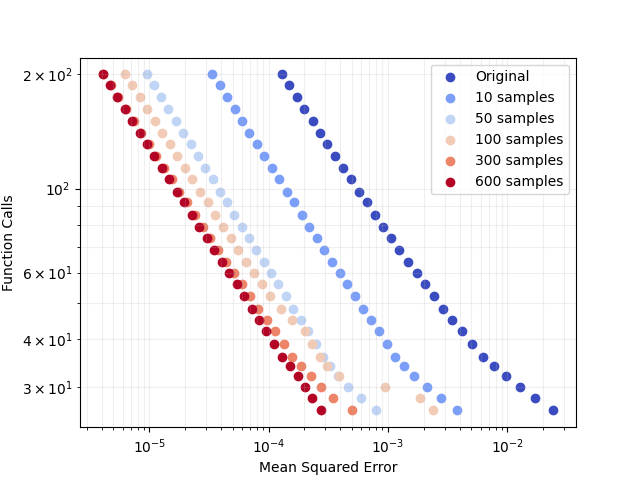}
\subcaption{Pseudo-invertible network}
\end{minipage}\hfill
\begin{minipage}{0.5\textwidth}
\centering
\includegraphics[width=\textwidth]{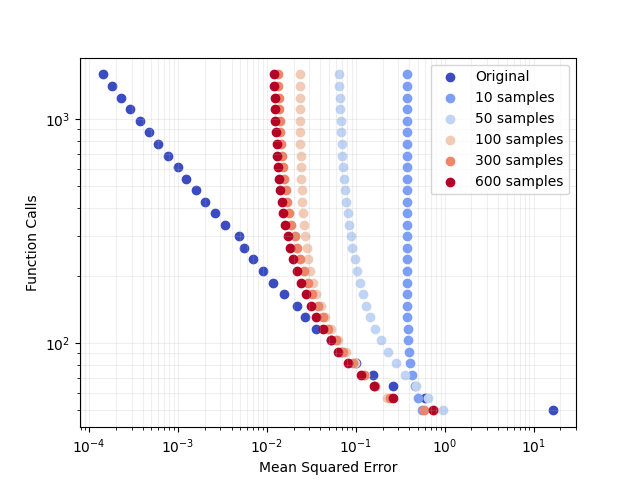}
\subcaption{PINODE}
\end{minipage}
    \caption{Euler solver function calls vs. final mean squared error for networks trained on various training sample sizes.}
    \label{fig:euler_sample_study}
\end{figure}

\textbf{Sample Size.} We trained networks to accelerate the linear system using training sets ranging from $10$ to $600$ uniformly distributed random points, then tested on trajectories outside the training range. Fig. \ref{fig:euler_sample_study}(a) shows that with $10$ sample points, Euler latent simulations were nearly $2$x faster than the original, and the speedup increased steadily until converging around $5$x with $300$ samples. Fig. \ref{fig:euler_sample_study}(b) shows the PINODE also achieves higher accuracy with more training samples but it reaches a barrier, and cannot achieve errors lower than $10^{-2}$ because the test trajectories were outside the range of training data provided.

\begin{figure}[hbt!]
\centering
\begin{minipage}{0.49\textwidth}
\centering
\includegraphics[width=\textwidth]{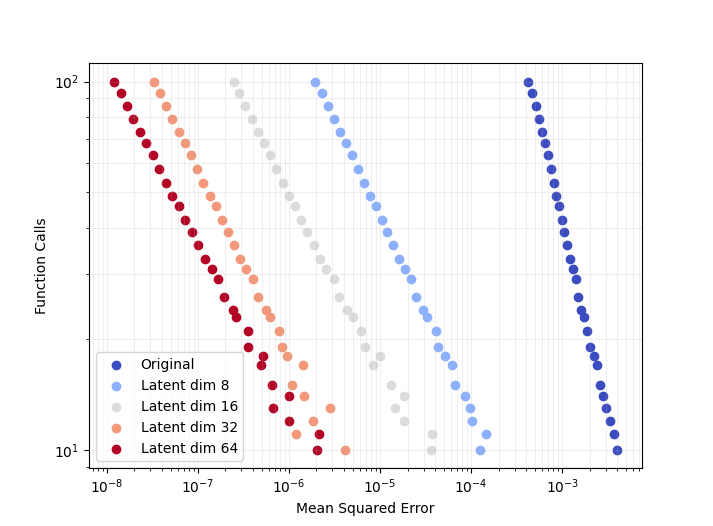}
\caption{Effect of latent dimension on Euler solver performance for vortex particle dynamics.}
\label{fig:euler_dim_study}
\end{minipage}\hfill
\begin{minipage}{0.49\textwidth}
\centering
\includegraphics[width=\textwidth]{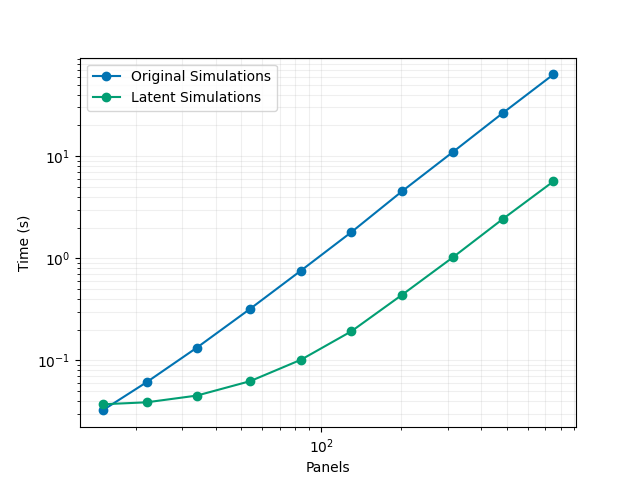}
\caption{Total computation time vs. VLM panel count for original and latent simulations.}
\label{fig:panel_count_study}
\end{minipage}
\end{figure}

\textbf{Latent Dimension.} We trained networks with latent dimensions of $8$, $16$, $32$, and $64$ for the vortex particle example, and observed how this parameter affected performance. Fig. \ref{fig:euler_dim_study} shows that even a latent dimension of $8$ (matching the original system size) results in over $10$x speedup for Euler simulations in this example. Speedups improve significantly at $16$ dimensions with diminishing returns beyond $32$. The higher the latent dimension, the better the network is able to speed up simulations. This is likely because more latent dimensions provide more representational flexibility for the network, which can learn to represent rapidly-changing signals as high-dimensional slow signals.

\textbf{Problem Complexity.} We simulated the VLM system with varying levels of panel discretization to measure how the network speedups scale with problem complexity. Fig. \ref{fig:panel_count_study} compares total simulation time versus VLM panel count for latent and original simulations of the same accuracy. The original simulations used an RK4 solver with a time-step of $0.06$s, and the latent simulations used an Euler solver with a time-step of $0.2$s.  These solver settings were chosen because both resulted in a final accuracy of about $2 \times 10^{-4}$, and remained at a similar level through the entire experiment. At low panel counts, speedups are limited, and latent simulations can be slower than the original. Speedups increase steadily with panel count, eventually converging to about $10$x. This demonstrates a fundamental tradeoff: the neural network adds computational overhead to each function call through the Jacobian-vector-product calculation, but when the original equations are already costly, this overhead becomes negligible and speedups are more pronounced.

\section{Conclusion}

In this work, we developed a method for speeding up the simulation of ordinary differential equations.  We used pseudo-invertible networks to map system states into a higher-dimensional latent space, where the governing equations were derived directly from the original system equations and the chain rule.  The network was trained to penalize the norm of the latent system's Jacobian matrix, which smoothed out the system dynamics and allowed ODE solvers to safely take larger steps. This method was applied to several example ODEs, and demonstrated $3$x to $20$x speedups over the original simualtions while maintaining the ability to generalize beyond training data with arbitrary accuracy. 


%
%

\paragraph{Limitations:} The main limitations of this approach are

\begin{enumerate}

    \item The method requires full knowledge of the system dynamics $\frac{d\mathbf{x}}{dt} = f(\mathbf{x})$ for the derivation of the latent equations.  For scenarios where the equations are partially or completely unknown, traditional data-driven methods are necessary. 
    \item The derived latent dyamic equations are more computationally complex than the original because they involve evaluating the original equations, as well as taking the derivative of the neural network.  Therefore it is well suited for problems where the original equations are expensive, and the overhead of the network is negligible.  For systems that are already fast to evaluate, the method will not perform as well. 
        
\end{enumerate}

Future work in this direction includes the application of this method to new physical systems, including complex high-dimensional systems like fluid flow or robotics.  There is also the possibility of extending the method to other forms of differential equations such as implicit differential equations, differential algebraic equations, and partial differential equations.  Ultimately, this work presents a framework for applying deep neural networks to build upon our best physical models of the world and achieve efficient, general solutions to complex differential equations.

%
%
%

\bibliographystyle{plain}
\bibliography{refs}

\appendix

\section{Network Architecture and Training Details}

This appendix provides complete specifications for the pseudo-invertible networks used in each experiment. All networks use the same basic architecture consisting of a pseudo-invertible linear layer followed by affine coupling layers.


\begin{table}[h]
\centering
\caption{Network architecture and training hyperparameters for the linear system example (Section 3.1).}
\label{tab:linear_system}
\begin{tabular}{ll}
\toprule
\textbf{Parameter} & \textbf{Value} \\
\midrule
Input dimension ($n$) & 3 \\
Latent dimension ($m$) & 64 \\
Number of coupling layers & 6 \\
Coupling layer height & 64 \\
Coupling layer depth & 3 hidden layers \\
Activation function & SiLU \\
Optimizer & Adam \\
Learning rate & 0.001 \\
Training epochs & 30,000 \\
Number of sample points ($N$) & 600 \\
Random JVP directions ($k$) & 8 \\
\bottomrule
\end{tabular}
\end{table}


\begin{table}[h]
\centering
\caption{Network architecture and training hyperparameters for the vortex particle system (Section 3.2).}
\label{tab:vortex_particles}
\begin{tabular}{ll}
\toprule
\textbf{Parameter} & \textbf{Value} \\
\midrule
Input dimension ($n$) & 8 \\
Latent dimension ($m$) & 64 \\
Number of coupling layers & 8 \\
Coupling layer hidden units & 64 \\
Coupling layer depth & 5 hidden layers \\
Activation function & SiLU \\
Optimizer & AdamiW \\
Learning rate & 0.001 \\
Training epochs & 50,000 \\
Number of training trajectories & 40 \\
Random JVP directions ($k$) & 4 \\
\bottomrule
\end{tabular}
\end{table}


\begin{table}[h]
\centering
\caption{Network architecture and training hyperparameters for the VLM aircraft model (Section 3.3).}
\label{tab:vlm}
\begin{tabular}{ll}
\toprule
\textbf{Parameter} & \textbf{Value} \\
\midrule
Input dimension ($n$) & 4 \\
Latent dimension ($m$) & 32 \\
Number of coupling layers & 8 \\
Coupling layer hidden units & 64 \\
Coupling layer depth & 3 hidden layers \\
Activation function & SiLU \\
Optimizer & AdamW \\
Learning rate & 0.0005 \\
Training epochs & 50,000 \\
Number of training trajectories & 20 \\
Random JVP directions ($k$) & 8 \\
\bottomrule
\end{tabular}
\end{table}

\newpage

\end{document}